\documentclass[a4paper, 11pt]{article}
\usepackage{amsthm}
\usepackage{amsfonts}
\usepackage{latexsym}
\usepackage{epic}
\usepackage{graphics}
\usepackage{graphicx}%
\usepackage{epsf, epstopdf}
\usepackage{etex}
\usepackage{pinlabel}

\usepackage[top=2.5cm,bottom=2.5cm,left=2.5cm,right=2.5cm]{geometry}
\usepackage{graphicx}
\usepackage{amsmath,amsthm,amssymb,lineno}

\usepackage{color, ifpdf}  
\usepackage[dvipsnames]{xcolor}            % Ö§³Ö²ÊÉ«
\usepackage{geometry}
\usepackage{float}
\usepackage{indentfirst}        % Ê×ÐÐËõœøºê°ü
\usepackage{amsmath,amssymb,bm} % AMS mathºê°üÓëÊýÑ§·ûºÅŒÓŽÖ
\usepackage{cases}              % ÊýÑ§¹«Êœ
\usepackage{setspace}
\usepackage{graphicx,booktabs,multirow}
\usepackage{enumerate,mathdots}

\usepackage{latexsym}
\usepackage{graphicx,color, booktabs,multirow}
\usepackage{tikz,chemfig}
\usepackage{graphicx,booktabs,multirow,mathrsfs}
\usepackage{appendix}
\usepackage{yhmath}
\usepackage{setspace}

\usetikzlibrary{backgrounds}
\usetikzlibrary{arrows}
\usetikzlibrary{shapes,shapes.geometric,shapes.misc}
\pgfdeclarelayer{edgelayer}
\pgfdeclarelayer{nodelayer}
\pgfsetlayers{background,edgelayer,nodelayer,main}
\tikzstyle{none}=[inner sep=0mm]
\tikzstyle{every loop}=[]

\tikzstyle{dotted}=[dash pattern=on \pgflinewidth off 2pt]
\tikzstyle{dashed}=[dash pattern=on 3pt off 3pt]

\newcommand \tikzp[2]
{
	\begin{center}
		\begin{tikzpicture}[scale=#1]
			#2
		\end{tikzpicture}
	\end{center}
}

\tikzstyle{new style 0}=[fill=black, draw=black, shape=circle]
\tikzstyle{red style 1}=[fill=red, draw=black, shape=circle]
\tikzstyle{blue style 2}=[fill=blue, draw=black, shape=circle]
\tikzstyle{white style 4}=[fill=white, draw=black, shape=circle]
\tikzstyle{bklack style 5}=[fill=black, draw=black, shape=rectangle]
\tikzstyle{red style 3}=[fill=red, draw=black, shape=rectangle]
\tikzstyle{yellow style 7}=[fill=yellow, draw=black, shape=rectangle]
\tikzstyle{new style 8}=[fill={rgb,255: red,0; green,132; blue,0}, draw={rgb,255: red,0; green,131; blue,0}, shape=circle]

\tikzstyle{new edge style 0}=[-]
\tikzstyle{new edge style 1}=[-, draw=red]
\tikzstyle{new edge style 2}=[-, draw=blue]
\tikzstyle{new edge style 3}=[-, draw={rgb,255: red,0; green,156; blue,0}]

\allowdisplaybreaks

\numberwithin{equation}{section}

\newcounter{countcase}

\newcounter{countclaim}
\def\inclaim{\addtocounter{countclaim}{1}
{\noindent {\bf Claim \thecountclaim}: }}

\def \proof {\noindent {\it Proof}. }
\newcommand{\proofend}
{
\setcounter{countclaim} {0}
\setcounter{countcase} {0}
{\hfill$\Box$}
}

\newcommand{\claimend}{{\hfill $\natural $}}

\def \T {{\mathcal {T}}}
\def \C {{\mathcal {C}}}
\def \S {\mathcal{S}}

\def \seq {{Seq_{-}or}}
\def \P {{\mathscr P}}

\newcommand \equ[2]
{
\begin{equation}\label{#1}
#2
\end{equation}
}

\newcommand \eqn[2]
{
\begin{eqnarray}\label{#1}
#2
\end{eqnarray}
}

\newcommand \sq[2]
{
	\S_{#2}(#1)
}

\newcommand \msq[2]
{
	\overline{\S_{#2}}(#1)
}

\newcommand \ellq[2]
{
\ell_{#2}(#1)
}

\begin{document}

 \tikzstyle{cblue}=[circle, draw, thin,fill=blue!20, scale=0.5]
%%%%%%===== ±êÌâÃû³ÆÖÐÎÄ»¯ =====
\renewcommand{\contentsname}{\LARGE \center Content}         %Ä¿\quad ÂŒ
                         %ÎÊÌâ

\newtheorem{theo}{Theorem}[section]
\newtheorem{lemm}[theo]{Lemma}
\newtheorem{deff}[theo]{Definition}
\newtheorem{coro}[theo]{Corollary}
\newtheorem{clm}{Claim}

\newcommand \Clm[2]
{
\begin{clm}\label{#1}
#2
\end{clm}
}

\title{
Anti-Ramsey numbers for trees in  
complete multi-partite graphs %and subtrees
}

\author{Meiqiao Zhang\thanks{Corresponding author. Email: nie21.zm@e.ntu.edu.sg and 
		meiqiaozhang95@163.com.} and Fengming Dong\thanks{Email: fengming.dong@nie.edu.sg (expired on 24/03/2027) and donggraph@163.com.}\\
\small National Institute of Education,
Nanyang Technological University, Singapore
}

\date{}
\maketitle

\begin{abstract}
Let $G$ be a complete multi-partite graph of order $n$.
In this paper, we consider the anti-Ramsey number 
$ar(G,\T_q)$ with respect to $G$ and 
the set $\T_q$ of trees with $q$ edges,
where $2\le q\le n-1$.
For the case 
$q=n-1$,
the result has been obtained by
Lu, Meier and Wang.
We will extend it to $q<n-1$.
We first show that $ar(G,\T_q)=\ell_q(G)+1$, 
where $\ell_q(G)$ is the maximum size of 
a disconnected spanning subgraph $H$ of $G$ with the property that 
any two components of $H$ together have at most $q$ vertices.
Using this equality, 
we obtain the exact values of $ar(G,\T_{q})$ 
for $n-3\le q\le n-1$.
We also compute $ar(G,\T_{q})$ by a simple algorithm 
when $(4n-2)/5\le q\le n-1$. 
\end{abstract}

%\tableofcontents

\section{Introduction
	\label{sec1}}

In this article, we consider simple graphs only. 
Given any graph $G$, let $V(G)$ and $E(G)$ denote the vertex set and edge set of $G$, and let $com(G)$ denote the number of its components.
If $G_1,G_2,\cdots, G_s$ are the components of $G$,
where $s=com(G)$,  
with $|V(G_1)|\ge \cdots \ge |V(G_s)|$, 
let $or_i(G)=|V(G_i)|$ for all  $i=1,2,\cdots, s$.  
Thus, $or_1(G)\ge or_2(G)\ge \cdots \ge or_s(G)$
and $or_1(G)+or_2(G)+\cdots + or_s(G)=|V(G)|$. 
Let $K_{p_1,p_2,\cdots,p_k}$ denote the complete 
$k$-partite graph whose partite sets' sizes 
are $p_1,p_2,\cdots,p_k$ respectively. 
For any $1\le r\le |V(G)|$, let 
$\P_r(G)$ be the family of $r$-element subsets of $V(G)$.
For any vertex $v\in V(G)$, 
let $E_G(v)$ be the set of edges in $G$ 
which are incident with $v$,
$N_G(v)$ be the set of vertices in $G$ which are 
adjacent to $v$ and 
$d_G(v)$ be the \textit{degree} of $v$ in $G$, i.e., the cardinality of $N_G(v)$. 
For any $S\subseteq V(G)$, 
let $E_G(S)=\bigcup_{v\in S}E_G(v)$
and let $G[S]$ be the subgraph of $G$ induced by $S$.

For a positive integer $t$, 
a \textit{$t$-edge-coloring} of $G$ is a surjective map from $E(G)$ to $\{1, 2, \cdots, t\}$. 
Note that an edge coloring here is actually a partition 
of $E(G)$, and it is probably not a proper edge coloring of $G$.
In an edge coloring of $G$, a subgraph $H$ of $G$ 
is called a \textit{rainbow} subgraph 
if the colors assigned to the edges in $H$ are pairwise distinct.

For a graph $G$ and a family $\C$ of graphs, 
the \textit{anti-Ramsey number} 
with respect to $G$ and $\C$,
denoted by $ar(G, \C)$,
is the maximum integer $t$ such that 
there is a $t$-edge-coloring of $G$ 
in which every rainbow subgraph 
is not isomorphic to any graph in $\C$. 
If no such edge coloring exists, define $ar(G, \C)$ to be zero.
The study of anti-Ramsey numbers was initiated by 
Erd\H{o}s, Simonovits and S\'os~\cite{4}. 
Since then, a lot of research papers 
on this topic have been published.
See~\cite{Fujita2014} for a survey and
~\cite{1, Fang2020, Xie2020, Yuan2021} for some recent development on specific $G$ or $\C$.

In particular, a common type of $\C$ is related to trees.
For example, the anti-Ramsey number for edge disjoint spanning trees has been studied thoroughly, the exact value of which has been obtained when $G$ is a complete graph~\cite{Jahan2016}, a complete bipartite graph~\cite{Jia2021}, and a complete multi-partite graph in~\cite{Lu2021} very recently. Moreover, the anti-Ramsey number for edge disjoint spanning trees in general graphs has also been determined in~\cite{Lu2021}.

From a different perspective, we will focus on another class of trees in the following.

This paper is motivated by the known results 
of $ar(G,\T_q)$ (see~\cite{3,5}),
where $G$ is a complete graph $K_n$ or 
a complete bipartite graph $K_{p_1,p_2}$ and
$\T_{q}$ is the set of subtrees in $G$ 
with exactly $q$ edges. %, where $2\le q\le n-1$. 
In~\cite{3}, Jiang and West determined and obtained the exact value of $ar(K_n,\T_q)$ by proving the equality that for $2\le q\le n-1$,  
\equ
{ne2-1}
{
ar(K_n,\T_q)=\ellq{K_n}{q}+1,
}
where
$\ellq{G}{q}$
is the maximum size of a disconnected spanning 
subgraph $H$ of $G$ in which every two components together 
have at most $q$ vertices (i.e., $or_1(H)+or_2(H)\le q$).
A result similar to (\ref{ne2-1}) for $ar(K_{p_1, p_2}, \T_q)$
was obtained by Jin and Li~\cite{5}, accompanying with exact values for certain cases. 

As a generalization, we consider 
the anti-Ramsey number 
$ar(G,\T_q)$, where $G$ is a complete multi-partite graph 
$K_{p_1,p_2,\cdots,p_k}$, $n=\sum_{i=1}^{k}p_{i}$ 
and $2\leq q\leq n-1$.
Note that the exact value of $ar(G,\T_{n-1})$ can be obtained as a corollary of the result in~\cite{Lu2021}.
%Recently,  Lu, Meier and Wang~\cite{Lu2021}
%determined $ar(G,\T_{n-1})$.
We will study the problem in a different approach.
We first extend the result (\ref{ne2-1}) to  
$ar(G,\T_q)$ for any complete multi-partite graph $G$ 
and $2\le q\le n-1$: 
%in Section~\ref{sec4}:
\equ{main1}
{
ar(G,\T_q)=\ellq{G}{q}+1,
}
transforming the study of $ar(G,\T_q)$ to that of $\ellq{G}{q}$, which is a seemingly more numerical invariant.

In this article, we calculate $ar(G,\T_q)$ 
via determining $\ellq{G}{q}$ 
for the two cases
$n-3\le q\le n-1$ and $(4n-2)/5\le q\le n-1$.
For both cases, %When $n-3\le q\le n-1$, 
%In order to determine $\ellq{G}{q}$ for a complete multi-partite graph $G$ and , 
we show that 
\equ{eq1-12}
{
\ellq{G}{q}=|E(G)|-
\min\limits_{S\in \P_{n-q+1}(G)}|E_G(S)|, 
}
%where $\P_r(G)$ is the family of $r$-element subsets of $V(G)$,
unless $(G,q)$ is one of the ordered pairs below
in the first case:
\equ{eq1-13}
{
(K_{3,3},4), (K_{4,3},4),  (K_{3,3,3},6).
}

The equality of (\ref{main1}) is established in 
Section~\ref{sec4}. 
Let $\sq{G}{q}$ be the set of 
disconnected spanning subgraphs $H$ of $G$
such that every two components of $H$ together have at most $q$ vertices and $\msq{G}{q}$ be the subset of graphs
in $\sq{G}{q}$ with
the maximum size (i.e., $\ellq{G}{q}$). 
In Section~\ref{sec06}, we show that 
%for a complete multi-partite graph $G$, 
$\msq{G}{q}$ contains a graph $H$ with $or_1(H)+or_2(H)=q$ and $or_2(H)=or_3(H)$.
In Section~\ref{sec6-1},
we consider the case 
$n-3\le q\le n-1$ and prove that 
(\ref{eq1-12}) holds,
unless $(G,q)$ is an ordered pair in (\ref{eq1-13}). 
In Section~\ref{sec6-2}, 
we get a conclusion that 
%For a complete multi-partite graph $G$, the value of 
$\min\limits_{S\in \P_{r}(G)}|E_G(S)|$
can be determined by a simple algorithm 
(i.e., Algorithm A in Section~\ref{sec6-2})
of repeatedly choosing vertices with the minimum degree.
Applying the results obtained in Sections~\ref{sec6-1}
and~\ref{sec6-2},
we explicitly express $ar(G,\T_q)$ 
for the case $n-3\le q\le n-1$
in Section~\ref{sec6-3}.
The case $(4n-2)/5\le q\le n-1$ is studied in the 
last section (i.e. Section~\ref{sec6-4}).
For this case, $\ellq{G}{q}$ can be determined 
by (\ref{eq1-12}) and thus it can be calculated by 
Algorithm A.

Let 
%$rank(G)$ denote the \textit{rank} of $G$, which is defined to be $rank(G)=|V(G)|-com(G)$, 
$\chi(G)$ denote the \textit{chromatic number} of graph $G$ and $\delta(G)$ be the \textit{minimum degree} among all vertices of $G$.
If $H$ is a proper subgraph of $G$, 
write $G\setminus H$ for $G[V(G)\setminus V(H)]$.
Especially, write $G\setminus x$ for $G[V(G)\setminus\{x\}]$.
Given a subset $E'$ of $E(G)$, 
let $G\setminus E'$ be the spanning subgraph of $G$ 
with edge set $E(G)\setminus E'$. 
If $E'=\{e'\}$, write $G\setminus e'$ for $G\setminus E'$.

\section{The anti-Ramsey number 
$ar(K_{p_1,\cdots,p_k},\T_q)$
	\label{sec4}
}

In this section,
we show that  
$ar(K_{p_1,\cdots,p_k},\T_q)
=\ellq{K_{p_1,\cdots,p_k}}{q}+1$.
We first introduce three lemmas.

\begin{lemm}\label{Crqlb}
	Let $G$ be any connected graph of order $n\ge 3$.  
	For any $q$ with $2\le q\le n-1$,
	\equ{ne2-3}
	{
		ar(G,\T_{q})\ge \ellq{G}{q}+1.
	}
\end{lemm}

\proof Let $H\in \msq{G}{q}$, i.e., $H\in \sq{G}{q}$ and
$|E(H)|=\ellq{G}{q}$.
Consider an $(\ellq{G}{q}+1)$-edge-coloring
of $G$ which assigns distinct colors to all edges in $H$ 
and a new color to 
all the edges in $E(G)\setminus E(H)$. We will show that it is an edge coloring without any rainbow subtree of $q$ edges.

Let $T$ be a subtree of $G$ with $q$ edges and $q+1$ vertices. 
By the definition of $\sq{G}{q}$, for any two components 
$H_1$ and $H_2$ of $H$, $|V(H_1)|+|V(H_2)|\le q$.
Then $T$ contains vertices from at least three components of $H$,
%$|V(T)\cap (V(G)\setminus V(H))|\ge 1$, 
thus $|E(T)\cap (E(G)\setminus E(H))|\ge 2$ and $T$ is not rainbow.

As there is no rainbow subtree of $q$ edges in this coloring,
the conclusion holds.
\proofend

Given a $t$-edge-coloring $c$ of a graph $G$, for any edge $e\in E(G)$, denote the color of $e$ by $c(e)$. A \textit{representing graph} of the $t$-edge-coloring $c$ is a rainbow spanning subgraph of $G$ with precisely $t$ edges.

\begin{lemm}\label{mainlem}
 Let $G$ be any connected graph and $c$ be a $t$-edge-coloring of $G$.
 Assume that $H$ is a representing graph of coloring $c$
 such that $H$ has a component $H_1$ with 
$|V(H_1)|\ge or_1(H')$ for every 
%the largest order among the components of all 
representing graph $H'$ of this coloring.
If $com(H)\ge 2$, % other than $H_1$, 
then $H_1$ has a bridge $b$. 
\end{lemm}

\proof Let $H_1, \cdots,H_s$ be 
all the components of $H$, where $s\ge 2$.

Since $G$ is connected, there exists an edge 
$e\in E(G)\setminus E(H)$
such that $e$ joins a vertex $u$ in $H_1$
to a vertex $v$ in $H_i$, where $2\le i \le s$.
As $e$ is not in $H$ and $H$ is a representing 
graph of coloring $c$, there must be a unique edge $b$ in $E(H)$
such that $c(b)=c(e)$. We will show $b$ is a bridge of $H_1$.

Let $H'$ be the spanning subgraph of $G$ obtained from $H$ by removing 
edge $b$ and adding edge $e$. 
$H'$ is also a representing graph for coloring $c$.
If either $b$ is not in $H_1$ or $b$ is in $H_1$
but $b$ is not a bridge of $H_1$, 
then $H'$ has a component including vertex $v$ and all the vertices of $H_1$, 
whose order is strictly larger than $|V(H_1)|$,
a contradiction to the selection of $H$. 

Hence $H_1$ has a bridge $b$.
\proofend

The following result due to Jiang and West~\cite{3}
will also be 
applied in the proof of Theorem~\ref{theoq}.

\begin{lemm}[\cite{3}]\label{lemmm}
	Every connected graph $G$ contains a vertex $w$ such that for 
each $e\in E(G)$, the component of $G-e$ containing $w$ has at least $|V(G)|/2$ vertices.
\end{lemm}

Now we can prove the main result in this section.
%$ar(K_{p_{1}, \cdots, p_{k}}, \T_{q})$ by $\ellq{K_{p_1,\cdots,p_k}}{q}$.

\begin{theo}\label{theoq}
For any positive integers $k, p_1,p_2,\cdots,p_k, q$ 
with $k\ge 2$ and $2\leq q\leq \sum_{i=1}^{k}p_{i}-1$, 
\equ{theoq-e1}
{
ar(K_{p_{1}, \cdots, p_{k}}, \T_{q})=\ellq{K_{p_1,\cdots,p_k}}{q}+1.
}
\end{theo}

\proof %\noindent \textit{Proof of Theorem~\ref{theoq}}:
Let $n=\sum_{i=1}^{k}p_{i}$ and $X_1, X_2,\cdots, X_k$ be the partite sets of  $K_{p_{1}, \cdots, p_{k}}$
with $|X_i|=p_i$ for all $i=1,2,\cdots,k$.
According to Lemma~\ref{Crqlb}, 
$ar(K_{p_{1}, \cdots, p_{k}}, \T_{q})\ge \ellq{K_{p_1,\cdots,p_k}}{q}+1.$
Thus, we only need to prove $ar(K_{p_{1}, \cdots, p_{k}}, \T_{q})\le \ellq{K_{p_1,\cdots,p_k}}{q}+1.$

Let $t=ar(K_{p_{1}, \cdots, p_{k}}, \T_{q})$ and $G$ denote $K_{p_1,\cdots,p_k}$. We will show that for any $t$-edge-coloring of $G$, if there is no rainbow subtree of $q$ edges, then $t\le \ellq{G}{q}+1$.

Now let $c$ be any $t$-edge-coloring 
of $G$ without rainbow trees of $q$ edges and choose a representing graph $H$ %$H=(X_1, \cdots, X_k)$
of coloring $c$ %with a component $H_1$ 
such that $or_1(H)\ge or_1(H')$
%$H$ has a component $H_1$ with the largestorder among the components of all 
for every representing graph $H'$ of this coloring.
%such that $H$ has a component $H_1$ with the largestorder among the components of all representing graphsof this coloring. 

Let $H_1$ be a component of $H$
with $|V(H_1)|=or_1(H)$.
By Lemma~\ref{Crqlb}, $|E(H)|=t\ge 1$, thus $|V(H_1)|\ge 2$.
Moreover, since there is no rainbow tree of $q$ edges in coloring $c$,
$|V(H_1)|\le q\leq n-1$, 
implying that $H$ has more than one component. %s other than $H_1$. 
Let $H_2,\cdots,H_s$ be the components of $H$ 
other than $H_1$, with $|V(H_i)|=or_i(H)$ for $i=2,\cdots,s$, where $s\ge 2$. 

%Then there is a vertex $v\in N_{H_1}(u)\cap X_{j}$ for some $1\leq j\leq k$, where $j\neq i$.
%such that $u$ and $v$ are adjacent in $H_1$. 

We will show $t\le \ellq{G}{q}+1$ by the following two cases.

\inclaim %$V(H_p)\subseteq X_i$ for all $p=2, 3,\cdots, s$.
If $|V(H_2)|=1$, then $t\le \ellq{G}{q}+1$.

% for all $p=2, 3,\cdots, s$.

%Clearly, $|V(H_p)|=1$ for all $p=2,3,\cdots, s$.

\proof
As $G$ is connected, 
there must be an edge $e=xz\in E(G)\setminus E(H)$, 
where $z\in V(H_1)$ and $x\in V(G)\setminus V(H_1)$.

\begin{figure}[!ht]
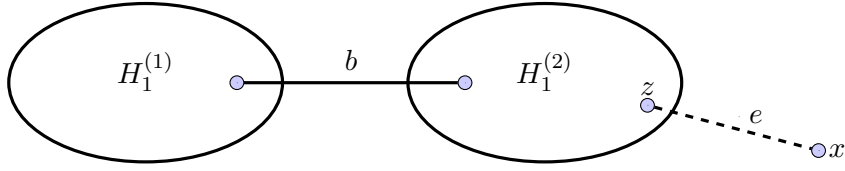

\tikzp{1.5}
{ 		\draw[black, very thick] (-4,0) ellipse (1.2 and 0.7);
		\draw[black, very thick] (-0.5,0) ellipse (1.2 and 0.7);

		\foreach \place/\y in {{(-3.2,0)/1}, {(-1.2,0)/2},{(0.4,-0.2)/3}, {(1.9,-0.6)/4}}
		\node[cblue] (b\y) at \place {};
		
		\filldraw[black] (b1) circle (0pt);
		\filldraw[black] (b2) circle (0pt);
		\filldraw[black] (b3) circle (0pt)node[anchor=south] {$z$};
        \filldraw[black] (b4) circle (0pt)node[anchor=west] {$x$};
        \filldraw[black] (1.2, -0.3) circle (0pt)node[anchor=west] {$e$};

		\draw[black, very thick] (b1) -- (b2);
		\draw[black, dashed, very thick] (b3) -- (b4);
		\node [style=none] (cap1) at (-4, 0.1) {$H_1^{(1)}$};
		\node [style=none] (cap2) at (-0.5, 0.1) {$H_1^{(2)}$};
		\node [style=none] (cap3) at (-2.2, 0.2) {$b$};
		
}
\caption{$H_1^{(1)}$ and $H_1^{(2)}$ are the components of 
		$H_1\setminus b$, and $x\in V(G)\setminus V(H_1)$}
\label{f1}
\end{figure}

By Lemma~\ref{mainlem},
there is a bridge $b\in E(H_1)$ with $c(e)=c(b)$.
Denote the two connected components of $H_1\setminus b$ by $H_1^{(1)}$ and $H_1^{(2)}$, as shown in Figure~\ref{f1}.

Let $H^*$ be the spanning subgraph of $G$ obtained from $H$ by deleting edge $b$. Note that the components of $H^*$ are 
$H_1^{(1)}, H_1^{(2)}, H_2,\cdots, H_s$.

Since $or_2(H)=|V(H_2)|=1$, 
we have 
$|V(H_p)|=1$ for all $p=2,3,\cdots, s$.
Thus, 
\equ{theoq-e2}
{
%|V(G_1^1)|+|V(K)|\leq|V(G_1^1)|+|V(K)|\leq
 or_1(H^*)+or_2(H^*)=|V(H_1^{(1)})|+|V(H_1^{(2)})|=|V(H_1)|\leq q.
} 
Hence $H^*\in \sq{G}{q}$ and 
%\equ{ne1-3}{
$t=|E(H)|=|E(H^*)|+1\le  \ellq{G}{q}+1$.
Claim 1 holds.
\claimend
%}

\inclaim If $|V(H_2)|\ge 2$, then $t\le \ellq{G}{q}+1$.
%There exists $w\in V(H_p)\cap X_{\alpha}$ for some $p$ and $\alpha$ with $2\le p\le s$, $1\le \alpha\le k$ and $\alpha\ne i$.

\proof
By Lemma~\ref{lemmm}, $H_1$ contains a vertex $u$, 
such that for each $e\in E(H_1)$, the component of 
$H\setminus e$ containing $u$ has at least $|V(H_1)|/2$ vertices. 
Assume $u\in V(H_1)\cap X_{i}$ for some $i$, where $1\leq i\leq k$.

Since $H_2$ is connected and $|V(H_2)|\ge 2$,  
we have $V(H_2)\not\subseteq X_i$.
Let $w\in V(H_2)\setminus X_i$.
Then $uw\in E(G)$.

\begin{figure}[!ht]
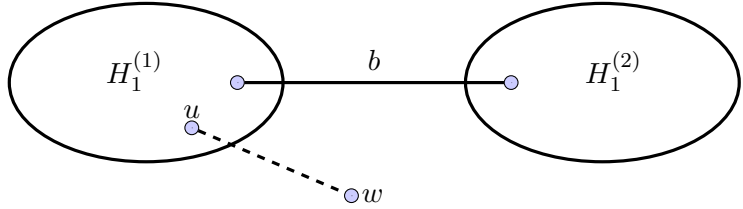

\tikzp{1.5}
	{ 	\draw[black, very thick] (-2,0) ellipse (1.2 and 0.7);
		\draw[black, very thick] (2,0) ellipse (1.2 and 0.7);

		\foreach \place/\y in {{(-1.2,0)/1}, {(1.2,0)/2},{(-1.6,-0.4)/3}, {(-0.2,-1)/4}}
		\node[cblue] (b\y) at \place {};
		
		\filldraw[black] (b1) circle (0pt);
		\filldraw[black] (b2) circle (0pt);
		\filldraw[black] (b3) circle (0pt)node[anchor=south] {$u$};
		\filldraw[black] (b4) circle (0pt)node[anchor=west] {$w$};
		
		\draw[black, very thick] (b1) -- (b2);
		\draw[black, dashed, very thick] (b3) -- (b4);
		\node [style=none] (cap1) at (-2.1, 0.1) {$H_1^{(1)}$};
		\node [style=none] (cap2) at (2.1, 0.1) {$H_1^{(2)}$};
		\node [style=none] (cap3) at (0, 0.2) {$b$};
		
}
\caption{$H_1^{(1)}$ and $H_1^{(2)}$ are the components of 
		$H_1\setminus b$, and $w\in V(H_2)$}
\label{f2}
\end{figure}

By Lemma~\ref{mainlem}, there is a bridge $b\in E(H_1)$ with $c(uw)=c(b)$. 
Denote the two connected components of $H_1\setminus b$ 
by $H_1^{(1)}$ and $H_1^{(2)}$ with $u \in H_1^{(1)}$,
as shown in Figure~\ref{f2}. %Assume $u \in H_1^{(1)}$.
%By Lemma~\ref{lemmm}, $|V(H_1^{(1)})|\ge |V(H_1^{(2)})|$.

Let $H'$ be the graph obtained from $H$ by deleting edge $b$ and adding edge $uw$. $H'$ is also a representing graph of coloring $c$
with a component which consists of vertices in both $H_1^{(1)}$ and $H_2$. 
By the assumptions of $H$ and $H_1$, 
we have 
$$
|V(H_1^{(1)})|+|V(H_2)|\leq |V(H_1)|=|V(H_1^{(1)})|+|V(H_1^{(2)})|,
$$
implying that 
$|V(H_1^{(2)})|\ge or_2(H)$. 
%$|V(H_1^{(2)})|\ge |V(H_2)|$. 

Let $H^*$ be the spanning subgraph of $G$ obtained from $H$
by deleting edge $b$. 
Note that $H^*$ has components 
$H_1^{(1)}, H_1^{(2)}, H_2,\cdots, H_s$. 

Due to the choice of $u$, $|V(H_1^{(1)})|\geq |V(H_1^{(2)})|$.
Hence $or_1(H^*)=|V(H_1^{(1)})|$ and $or_2(H^*)=|V(H_1^{(2)})|$.
%Hence $|V(H_1^{(1)})|\ge |V(H_1^{(2)})|\ge |V(H_2)|\ge \cdots \ge |V(H_{s})|$.
Together with the fact that 
$|V(H_1^{(1)})|+|V(H_1^{(2)})|=|V(H_1)|\le q$,  
we conclude that $H^*\in \sq{G}{q}$.
Hence $t=|E(H)|=|E(H^*)|+1\le \ellq{G}{q}+1$ and Claim 2 holds.
\claimend

By Claims 1 and 2, $ar(K_{p_{1}, \cdots, p_{k}}, \T_{q})\le \ellq{K_{p_1,\cdots,p_k}}{q}+1$,
and thus the result holds.
\proofend

\section{Preparation work \label{sec06}}

To get the exact value of $\ellq{K_{p_1,\cdots,p_k}}{q}$, we shall give some lemmas based on the properties of multi-partite graphs and $\ellq{G}{q}$ in this section.

\begin{lemm}\label{le6-4}
	Let $G$ be a complete multi-partite graph,
	and let $V_1$ and $V_2$ be disjoint subsets of $V(G)$.
	If $|V_1|\ge |V_2|$ and $V_1\cup V_2$ is not an independent set of $G$,
	then there exists a vertex $u\in V_2$ 
	such that 
	$|N_G(u)\cap V_1|\ge |N_G(u)\cap V_2|$,
	where the inequality is strict when either $|V_1|> |V_2|$ or 
	$\chi(G[V_2])<\chi(G[V_1\cup V_2])$.
\end{lemm}

\proof Since $G$ is a complete multi-partite graph,
$G':=G[V_1\cup V_2]$ is also a complete multi-partite graph.
Assume that $X_1,X_2,\cdots, X_k$ are the partite sets of $G'$,
where $k=\chi(G')$, 
and $Y_i=V_1\cap X_i$ and $Z_i=V_2\cap X_i$
for all $i=1,2,\cdots,k$.
Since $V_1\cup V_2$ is not an independent set of ${G'}$, $k\ge 2$.

Let $r:=\chi({G'}[V_2])$.
We may assume that 
$Z_i\ne \emptyset$ for all $i=1,2,\cdots,r$.
Then $Z_i=\emptyset$ and $Y_i=X_i$ for all $r+1\le i\le k$.

If $r=1$, then for each $u\in V_2$,
$|N_{G'}(u)\cap V_2|=0$ 
and $|N_{G'}(u)\cap V_1|=|V_1\setminus X_1|>0$, 
and the result is trivial.  Now assume that $r\ge 2$.
Note that 
\begin{eqnarray}\label{le6-4-e1}
	\frac 1{r-1} \sum_{i=1}^r (|V_2\setminus Z_i|)
	&=&\sum_{i=1}^r |Z_i|=|V_2|
	\nonumber \\
	&\le & |V_1| \nonumber \\
	&=&\sum_{i=1}^r |Y_i|+\sum_{i=r+1}^k |Y_i| \nonumber \\
	&=&\frac 1{r-1} \sum_{i=1}^r (|V_1\setminus Y_i|)
	-\frac 1{r-1}\sum_{i=r+1}^k |Y_i|\nonumber \\
	&\le &\frac 1{r-1} \sum_{i=1}^r (|V_1\setminus Y_i|),
\end{eqnarray}
where the inequality above is strict when either $|V_1|>|V_2|$ 
or $r<k$. Thus, there exists $i:1\le i\le r$ such that 
$$
|V_1\setminus Y_i|\ge |V_2\setminus Z_i|,
$$
where the inequality is strict when either $|V_1|>|V_2|$ 
or $r<k$.
Let $u\in Z_i$. Observe that 
\equ{}
{
	|N_{G'}(u)\cap V_1|-|N_{G'}(u)\cap V_2|=
	|V_1\setminus Y_i|-
	|V_2\setminus Z_i|\ge 0,
}
where the inequality is strict when either $|V_1|>|V_2|$ 
or $r<k$.

The result holds.
\proofend

\begin{lemm}\label{le6-40}
	Let $G$ be a complete multi-partite graph.
	For disjoint subsets $V_1$ and $V_2$ of $V(G)$, 
	if $\chi(G[V_2])=r$ and $\chi(G[V_1\cup V_2])=k$,
	then there exists  $u\in V_2$ 
	such that $|N_G(u)\cap V_1|\ge \frac{(r-1)|V_1|+k-r}{r}\ge \frac{(r-1)|V_1|}{r}$.
\end{lemm}

\proof Note that $G':=G[V_1\cup V_2]$ is a complete $k$-partite graph.
Assume that $X_1,X_2,\cdots, X_k$ are the partite sets of $G'$.
Let $Y_i:=V_1\cap X_i$ and $Z_i:=V_2\cap X_i$ for all $i=1,2,\cdots,k$.
Since $\chi(G[V_2])=r$, we may 
assume that $Z_i\ne \emptyset$  for all $i=1,2,\cdots,r$.
Then $Z_i=\emptyset$ and $Y_i=X_i$ for $r+1\le i\le k$.
Note that 
$$
|V_1|=\sum_{i=1}^k |Y_i|\ge (k-r)+\sum_{i=1}^r |Y_i|, 
$$
by which $|Y_i|\le \frac {|V_1|-k+r}r$ holds
for some $i$ with $1\le i\le r$.
Let $u\in Z_i$. 
As $G'$ is a complete $k$-partite graph with partite sets 
$X_1,X_2,\cdots, X_k$ and $u\in X_i$, 
$u$ is adjacent to all vertices in $V_1\setminus Y_i$,
implying that 
$$
|N_G(u)\cap V_1|=|V_1|-|Y_i|\ge |V_1|-\frac {|V_1|-k+r}r.
$$
The result holds.
\proofend

\begin{lemm}\label{le6-3}
	For any graph $Q$ with components $Q_1,Q_2,\cdots, Q_s$, 
	if $\chi(Q)=t$, then there exists 
	one vertex $u_i$ in $Q_i$ for each $i\in [s]$ such that 
	\equ{le6-3-e1}
	{
		\sum_{i=1}^s d_Q(u_i)\le \frac {t-1}{t}|V(Q)|.
	}
\end{lemm}

\proof We first show it for $s=1$.
%that (\ref{le6-3}) holds when $Q$ is connected.
Let $U_1,U_2,\cdots,U_t$ be the color classes of a proper
$t$-coloring with $|U_1|\ge |U_2|\ge \cdots \ge |U_t|$.
Then, for each vertex $u\in U_1$, we have 
\equ{le6-3-e2}
{
	d_Q(u)\le \sum_{i=2}^{t} |U_i|\le \frac {t-1}{t}|V(Q)|.
}
%Now consider the case that $Q$ has non-trivial components $Q_1,Q_2,\cdots, Q_s$.
By (\ref{le6-3-e2}), there exists vertex $u_i$ in $Q_i$ 
for all $i\in [s]$ such that 
\equ{le6-3-e3}
{
	\sum_{i=1}^sd_{Q}(u_i)=
	\sum_{i=1}^sd_{Q_i}(u_i)
	\le \sum_{i=1}^s 
	\frac {\chi(Q_i)-1}{\chi(Q_i)}|V(Q_i)|
	\le \sum_{i=1}^s 
	\frac {t-1}{t}|V(Q_i)|
	\le \frac {t-1}{t}|V(Q)|.
}
\proofend

\begin{lemm}\label{le2-1}
	Let $G$ be a connected graph of order $n\ge 3$ 
	and $q$ be an integer with $2\le q\le n-1$.
	%\begin{enumerate} \item  
	For any $H\in \sq{G}{q}$, $com(H)\ge 3$ holds, 
	% with components $H_1,\cdots,H_s$, $s\ge 3$ holds, 
	and if $H\in \msq{G}{q}$,  %$\ellq{G}{q}=|E(H)|$,
	then 
	\begin{enumerate} 
		\item  %If $\ellq{G}{q}=|E(H)|$, then 
		each component of $H$ is a vertex-induced subgraph of $H$; and 
		%\item $rank(H)\le \min\{n-3,n-2n/q\}$.
		\item  when %$|V(H_1)|\ge |V(H_2)|\ge \cdots \ge |V(H_s)|$ 
		$or_1(H)\ge or_1(H')$ for every graph $H'\in \msq{G}{q}$,
		then either $or_1(H)+or_3(H)=q$ or 
		$|N_G(u)\cap V(H_2)|> |N_G(u)\cap V(H_1)|$
		for each $u\in V(H_2)$, where $H_1$ and $H_2$ 
		are components of $H$ with $|V(H_i)|=or_i(H)$ for $i=1,2$.
	\end{enumerate}
\end{lemm}

\proof 
Let $s=com(H)$ and $H_1,H_2,\cdots,H_s$ be the components of $H$ with $|V(H_i)|=or_i(H)$ for $i=1,2,\cdots,s$.
As $H\in \sq{G}{q}$, $|V(H_1)|+|V(H_2)|\le q<n$, implying that $s\ge 3$.

Now assume that $\ellq{G}{q}=|E(H)|$.
(i) is trivial. 

(ii). Assume $|V(H_1)|\ge or_1(H')$ for every $H'\in \msq{G}{q}$.
Suppose (ii) fails. 
Then  %$|V(H_2)|>|V(H_3)|$ 
$|V(H_1)|+|V(H_3)|<q$
and 
$|N_G(u)\cap V(H_1)|\ge |N_G(u)\cap V(H_2)|$
for some $u\in V(H_2)$.
 
Note that $|N_G(u)\cap V(H_2)|\ge 1$, i.e., $|V(H_2)|\ge 2$.
Otherwise, $|V(H_2)|=|V(H_3)|=1$, $|V(H_1)|\le q-2$ and $|E(H)|=|E(H_1)|$.  
%Since $n\le (q-2)+s-1\le (n-1-2)+s-1$, $s\ge 4$. 
Since $G$ is connected, there exists an edge $uv\in E(G)$, where $u\in V(H_1)$ and $v\in V(G)\setminus V(H_1)$. Let $H'$ be a spanning subgraph of $G$ with edge set $E(G[V(H_1)\cup \{v\}])$. Then $or_1(H')\le q-1$ and $or_2(H')=1$, thus $H'\in \sq{G}{q}$, a contradiction to the assumption of $H\in \msq{G}{q}$ as $|E(H')|>|E(H)|$.

Let $H'_1$ denote the subgraph $G[V(H_1)\cup \{u\}]$
%let $H'_2$ denote the graph $H_2\setminus \{u\}$,
and %let $H'$ be the graph obtained from $H$ by replacing $H_1$ and $H_2$ by $H'_1$ and $H'_2$, i.e., 
$H'$ be the graph obtained from $H$ 
by adding all edges in $\{uv: v\in N_G(u)\cap V(H_1)\}$
and deleting all edges in $\{uv: v\in N_G(u)\cap V(H_2)\}$.
Since $|N_G(u)\cap V(H_1)|\ge |N_G(u)\cap V(H_2)|\ge 1$, $H'_1$ is connected and $or_1(H')=|V(H'_1)|=|V(H_1)|+1$. Moreover,
$$
|E(H')|-|E(H)|=|N_G(u)\cap V(H_1)|
-|N_G(u)\cap V(H_2)|\ge 0,
$$
%and $H'$ has at least $s$ components, and it has exactly $s$ components if and only if $H'_2$ is connected. Furthermore, 
 and 
$or_2(H')\le \max\{|V(H_2)|-1,|V(H_3)|\}$.
Thus,  
$$
or_1(H')+or_2(H')\le 
\max\{|V(H_1)|+|V(H_2)|,|V(H_1)|+1+|V(H_3)|\}\le q,
$$
implying that  $H'\in \sq{G}{q}$.
As $|E(H')|\ge |E(H)|=\ellq{G}{q}$, 
we have $H'\in \msq{G}{q}$, a contradiction 
to the assumption of $H$ as $or_1(H')>or_1(H)$.
%However, $|E(H')|>|E(H)|=\ellq{G}{q}$, a contradiction.

(ii) holds.
\proofend

%For any two infinite sequences $(a_1,a_2,\cdots)$ and $(b_1,b_2,\cdots)$,write $(a_1,a_2,\cdots) \succeq (b_1,b_2,\cdots)$if $a_i<b_i$ always implies that $i>1$ and $a_j>b_j$ holds for some $j:1\le j<i$.
For any two finite sequences  $(a_1,a_2,\cdots,a_s)$ and $(b_1,b_2,\cdots,b_t)$,
write $(a_1,a_2,\cdots,a_s) \succeq (b_1,b_2,\cdots,b_t)$
if either $s=t$ and $a_i=b_i$ for all $i=1,2,\cdots,s$,
or there exists $i:1\le i\le \min\{s,t\}$ such that 
$a_i>b_i$ and $a_j=b_j$ for all $1\le j<i$.
If $a_1+a_2+\cdots+a_s=b_1+b_2+\cdots+b_t$ and 
all $a_i$'s and $b_j$'s are positive,   
then either $(a_1,a_2,\cdots,a_s) \succeq (b_1,b_2,\cdots,b_t)$
or $(b_1,b_2,\cdots,b_t) \succeq  (a_1,a_2,\cdots,a_s)$.

For any graph $H$, let $\seq(H)$ denote the sequence 
$(or_1(H), or_2(H),\cdots, or_s(H))$, where $s=com(H)$. 
If $H$ and $H'$ are spanning subgraphs of $G$, 
then either $\seq(H)\succeq \seq(H')$ or $\seq(H')\succeq \seq(H)$.
Obviously, if $\seq(H)\succeq \seq(H')$, then $or_1(H)\ge or_1(H')$.

\iffalse
For spanning subgraphs $H$ and $H'$ of a graph $G$,
write $H\succeq_{or} H'$ if 
$(or_1(H),\cdots, or_s(H)) \succeq (or_1(H'),\cdots, or_r(H'))$,
where $s=com(H)$ and $t=com(H')$.
\fi

%$(a_1,a_2,\cdots,a_s,0,0,\cdots) \succeq (b_1,b_2,\cdots,b_t,0,0,\cdots)$. 

\begin{theo}\label{co6-1}
Let $G$ be a complete multi-partite graph of order $n$ with at least two partite sets.
	Assume that $2\le q\le n-1$ and $H$ is a graph in $\msq{G}{q}$ 
	such that $or_1(H)\ge or_1(H')$ for every $H'\in \msq{G}{q}$. The following hold:
	\begin{enumerate}
		\item $or_1(H)+or_3(H)=q$; and
		\item if $\seq(H)\succeq \seq(H')$ holds for every $H'\in \msq{G}{q}$,
		 then $\seq(H)$ is a sequence of the following form:
		 $$
		 (h_1,\underbrace{h_2,\cdots, h_2}_t, h_3)
		 $$
		 where $t\ge 1, h_1\ge h_2\ge h_3$, and $h_3=h_2$ when $t=1$; or
		 a sequence of the following form:
		 $$
		 (h_1,\underbrace{h_2,\cdots,h_2}_{t_1},\underbrace{1,\cdots, 1}_{t_2})
		 $$
		 where $t_1,t_2\ge 2$ and $h_1\ge h_2\ge 2$.
	\end{enumerate} 
\end{theo}

\proof 
Let $s=com(H)$ and $H_1,H_2,\cdots,H_s$ be the components of $H$ with $|V(H_i)|=or_i(H)$ for $i=1,2,\cdots,s$, where $s\ge 3$ by Lemma~\ref{le2-1}.

(i). The result is trivial when $q=2$. Assume $q\ge 3$ in the following.

Suppose that $or_1(H)+or_3(H)<q$.
Then, %by the assumption of $H$, 
Lemma~\ref{le2-1} (ii) implies that 
$|N_G(u)\cap V(H_2)|>|N_G(u)\cap V(H_1)|$ for each $u\in V(H_2)$.
But, as $|V(H_1)|\ge |V(H_2)|$, 
due to Lemma~\ref{le6-4},
$V(H_1)\cup V(H_2)$ is an independent set in $G$. Thus $|V(H_1)|=|V(H_2)|=1$ and for $u\in V(H_2)$, $|N_G(u)\cap V(H_2)|=|N_G(u)\cap V(H_1)|=0$, a contradiction. (i) holds.

(ii). We first prove the following claim.

\inclaim for any $5\le b\le s$, if $or_{b-1}(H)<or_2(H)$, then
$or_{b-1}(H)=1$.  

\proof Suppose the claim fails. Then, there exists $b$ 
with $5\le b\le s$ such that $2\le or_{b-1}(H)<or_2(H)$,
i.e., $2\le |V(H_{b-1})| <|V(H_2)|$.

By Lemma~\ref{le6-4}, 
there exists $u\in V(H_b)$ such that 
$|N_G(u)\cap V(H_b)|\le |N_G(u)\cap V(H_{b-1})|$.
Let $H'$ be the graph obtained from $H$ 
by adding all edges in $\{uv: v\in N_G(u)\cap V(H_{b-1})\}$
and deleting all edges in $\{uv: v\in N_G(u)\cap V(H_{b})\}$.
%As $|N_G(u)\cap V(H_{b-1})|\le |N_G(u)\cap V(H_p)|$, 
Obviously, $|E(H')|\ge |E(H)|$.

Let $H'_{b-1}$ be the subgraph $G[V(H_{b-1})\cup \{u\}]$.
%and edge set $E(H_{b-1})\cup \{uv: v\in N_G(u)\cap V(H_{b-1})\}$.
Since $|V(H_{b-1})|\ge 2$ and $G$ is a complete $k$-partite graph, 
$H'_{b-1}$ is connected. 

Thus, the components of $H'$ are 
$H_1, H_2, \cdots, H_{b-2}, H'_{b-1}, H_{b+1},\cdots, H_s$
together with components of $H_b\setminus \{u\}$,
implying that $\seq(H)\not \succeq \seq(H')$. % and $\seq(H')\ne \seq(H)$.

By the given condition, $|V(H'_{b-1})|=|V(H_b)|+1\le |V(H_2)|$.
Thus, $H'\in \sq{G}{q}$.
As $|E(H')|\ge |E(H)|$, $H\in \msq{G}{q}$.
However, $\seq(H)\not \succeq \seq(H')$,  
a contradiction to the assumption of $H$.

Hence Claim 1 holds.
\claimend

Let $or_i(H)=h_i$ for $i=1,2$.
By the result in (i), we have $or_3(H)=h_2$.
If $or_{s-1}(H)=h_2$, 
then $\seq(H)=(h_1,\underbrace{h_2,\cdots, h_2}_{t},h_3)$,
where $h_3=or_s(H)\le h_2$.
If $or_{s-1}(H)<h_2$, then $s\ge 5$ and 
there exists $5\le b\le s$ such that
$or_{b}(H)\le or_{b-1}(H)<or_2(H)=h_2$ and $or_{b-2}(H)=or_2(H)$.
In this case, due to Claim 1, $or_{b-1}(H)=1$,
implying that 
$\seq(H)=(h_1,\underbrace{h_2,\cdots, h_2}_{t_1},\underbrace{1,1,\cdots,1}_{t_2})$,
where $t_1,t_2\ge 2$.

Hence (ii) holds.
\proofend

%By applying Corollary~\ref{co6-1},
For a complete multi-partite graph $G$ of order $n$ 
and $2\le q\le n-1$, 
let $H$ be a member in $\msq{G}{q}$ such that $\seq(H)\succeq \seq(H')$ 
holds for every $H'\in \msq{G}{q}$. 
By Theorem~\ref{co6-1}, 
$h_1+h_2=q$ and $h_1+2h_2\le n$, thus $1\le h_2\le \min\{q/2,n-q\}$ and we can express all the possible sequences according to the value of $h_2$.
For example, if $n=13$ and $q=8$, 
then, $1\le h_2\le \min\{8/2,13-8\}=4$, and 
$\seq(H)$ is one of the sequences below:
$$
(4,4,4,1), (5,3,3,2),(5,3,3,1,1), (6,2,2,2,1), (6,2,2,1,1,1), (7,1,1,1,1,1,1).
$$

\section{To find $\ellq{K_{p_1,\cdots,p_k}}{q}$
for $n-3\le q\le n-1$
\label{sec6-1}
}

In this section, 
we consider the case $n-3\le q \le n-1$,
and show that (\ref{eq1-12}) holds 
unless $(G,q)$ is an ordered pair in (\ref{eq1-13}).
%(i.e., Theorem~\ref{th6-2}).}
%$\ellq{K_{p_1,\cdots,p_k}}{q}$ 
%we shall give a lower bound for $\ellq{K_{p_1,\cdots,p_k}}{q}$ in (\ref{eq6-8}) and determine $\ellq{K_{p_1,\cdots,p_k}}{q}$ when $n-3\le q \le n-1$ in Theorem~\ref{th6-2}.

Let $G$ be a connected graph of order $n$ and $q$ be an integer with 
$2\le q\le n-1$. 
%For any $0\le r\le n$, let $\P_r(G)$ be the family of $r$-element subsets of $V(G)$.
Recall that for any $S\subseteq V(G)$, 
$E_G(S)=\bigcup_{v\in S}E_G(v)$.
Clearly, 
\equ{eq6-08}
{
	|E(G)|-\min_{S_0\in \P_{r}(G)}
	|E_G(S_0)|
	%\left | \bigcup_{v\in S}E_G(v)\right |
	=\max_{S\in \P_{n-r}(G)}|E(G[S])|.
}
For any $S\in \P_{q-1}(G)$,
the spanning subgraph $H$ of $G$ with edge set 
$E(G[S])$ has at least
 $n-q+2$ components,
i.e., $G[S]$ and $n-q+1$ trivial components,  
and thus it belongs to $\sq{G}{q}$.
Hence 
\equ{eq6-8}
{
	\ellq{G}{q}
	\ge \max_{S\in \P_{q-1}(G)}|E(G[S])|
	=|E(G)|-
	\min_{S_0\in \P_{n-q+1}(G)}
	|E_G(S_0)|
	%\left | \bigcup_{v\in S_0}E_G(v)\right |.
}
It can be verified that the inequality of (\ref{eq6-8}) is strict 
when $(G, q)$ is one of the following ordered pairs:
\equ{eq6-13}
{
	(K_{3,3},4), (K_{4,3},4),  (K_{3,3,3},6).
}
For example, if $G=K_{3,3,3}$ and $q=6$, 
then $\ell_q(G)\ge 9$, as $\sq{G}{q}$ has a graph consisting of 
three components each of which is isomorphic to $K_3$,
while each $5$-vertex subgraph of $G$ has at most $8$ edges.
Similarly, if $G=K_{3,3}$ or $K_{4,3}$ and $q=4$,
each $3$-vertex subgraph of $G$ has at most $2$ edges,
while  $\sq{G}{q}$ has a graph containing $3$ edges.
%$|E(G[S])|\le 8$ for each $S\in \P_5(G)$. 

Now we shall show that
if $G$ is a complete multi-partite graph of order $n$
and $\max\{2,n-3\}\le q\le n-1$, then 
the equality of (\ref{eq6-8}) holds if and only if 
$(G, q)$ is not an ordered pair in (\ref{eq6-13}).
Thus, by applying Theorem~\ref{theoq}, 
$ar(G,\T_q)$ can be determined.

\begin{theo}\label{th6-2}
Let $G$ denote $K_{p_1,p_2,\cdots,p_k}$
and $n=p_1+p_2+\cdots+p_k$,
where $k\ge 2$ and $p_1\ge p_2\ge \cdots\ge p_k\ge 1$.
For $n-3\le q\le n-1$, %$1\le r\le 3$,
if $q\ge 2$ and $(G,q)$ is not an ordered pair in (\ref{eq6-13}),
then
\equ{eq6-9}
{
\ell_q(G)= |E(G)|
-\min_{S\in \P_{n-q+1}(G)}
|E_G(S)|.
%\left |  \bigcup_{v\in S}E_G(v)\right |.
}
\end{theo}

\proof Suppose that (\ref{eq6-9}) is not true. 
Then, the inequality of (\ref{eq6-8}) is strict. 

Assume that $(G,q)$ is not one of the ordered pairs in (\ref{eq6-13}) and $H$ is a graph in $\msq{G}{q}$ 
such that $or_1(H)\ge or_1(H')$ for every $H'\in \msq{G}{q}$. 
Let $H_1, H_2,\cdots, H_s$ be the components of $H$
with $|V(H_i)|=or_i(H)$ for $i=1,2,\cdots,s$, where $s\ge 3$ by Lemma~\ref{le2-1}.

Let $r=n-q$. Then $1\le r\le 3$.
%By Lemma~\ref{le2-1},
%$q-2\le rank(H)\le n-3$.As $s=n-rank(H)$, 
%$3\le s\le n-q+2=r+2$.

\inclaim $|E(H)|>|E(H_0)|$ holds for every 
subgraph $H_0$ of $G$ with $|V(H_0)|\le q-1$ ($=n-r-1$).

\proof
As $|E(H)|=\ell_q(G)$ and (\ref{eq6-9}) fails, 
the inequality of (\ref{eq6-8}) is strict,
and thus the claim follows.
\claimend

\inclaim $|V(H_2)|\ge 2$.

\proof
Suppose that $|V(H_2)|=1$. 
Observe that $|E(H_1)|=|E(H)|$.
As $H\in \sq{G}{q}$, we have $|V(H_1)|\le q-|V(H_2)|=q-1$.
It contradicts Claim 1.
Claim~\thecountclaim\ holds. 
\claimend

\inclaim $|V(H_1)|+|V(H_2)|=q$
and
$|V(H_2)|=|V(H_3)|\le r$.

\proof
Claim~\thecountclaim\ follows from Theorem~\ref{co6-1} (i).
\claimend

\inclaim $r>1$.

\proof
If $r=1$, then by Claim 3, $|V(H_2)|=1$,
a contradiction to Claim 2, thus Claim~\thecountclaim\ holds.
\claimend 

\inclaim 
$r\ne 2$.

\proof
Suppose that $r=2$.
By Claims 2 and 3, $|V(H_2)|=|V(H_3)|=2$ and $|V(H_1)|+|V(H_2)|+|V(H_3)|=n$, thus $s=3$, both $H_2$, $H_3$ are isomorphic to $K_2$ and $|E(H)|=|E(H_1)|+2$.

We are now going to show that $|V(H_1)|=2$. 
Suppose that $|V(H_1)|>|V(H_2)|$.
As $|V(H_2)|=2$, by Lemma~\ref{le6-4}, 
there exists a vertex $u\in V(H_2)$ such that 
$|N_G(u)\cap V(H_1)|\ge 1+|N_G(u)\cap V(H_2)|=2$.

Let $H_0$ be the subgraph $G[V(H_1)\cup \{u\}]$.
Observe that $|V(H_0)|=|V(H_1)|+1=n-3=n-r-1$ and 
$$
|E(H_0)|\ge |E(H_1)| +2=|E(H)|,
$$
a contradiction to Claim 1.

Hence $|V(H_1)|=2$, implying that $|E(H)|=3$ and $n=|V(H)|=6$.
Thus, $n-r-1=3$.
If $k\ge 3$, then $G$ contains a subgraph $H_0$ isomorphic to $K_3$,
a contradiction to Claim 1. 
Thus $k=2$. The three edges in $H$ 
form a perfect matching of $G$, 
implying that $G\cong K_{3,3}$,
a contradiction to the assumption of $G$.

Claim~\thecountclaim\ holds.
\claimend

\inclaim If $r=3$, then $s\le 3$.

\proof
Suppose that $r=3$ and $s\ge 4$. 
By Claims 2 and 3, $2\le |V(H_2)|=|V(H_3)|\le 3$ and $n-|V(H_1)|-|V(H_2)|=r=3$, implying $s\le 4$. Thus $s=4$,
$|V(H_2)|=|V(H_3)|=2$ and $|V(H_4)|=1$.

We are now going to show that $|V(H_1)|=2$. 
Suppose that $|V(H_1)|\ge 3>|V(H_2)|$.

As $|V(H_2)|=2$, by Lemma~\ref{le6-4}, 
there exists a vertex $u\in V(H_2)$ such that 
$|N_G(u)\cap V(H_1)|\ge 1+ |N_G(u)\cap V(H_2)|=2$.
Let $H_0$ be the subgraph $G[V(H_1)\cup \{u\}]$.
Observe that $|V(H_0)|=|V(H_1)|+1=n-4=n-r-1$ and
$$
|E(H_0)|\ge |E(H_1)|+2=|E(H)|,
$$
a contradiction to Claim 1.

Hence $|V(H_1)|=2$, implying that 
$|E(H)|=3$ and $n=7$.

Suppose that $k\ge 3$.  
Let $H_0$ be a subgraph isomorphic to $K_3$.
Observe that $|E(H_0)|= 3=|E(H)|$ and
$
|V(H_0)|=3=n-4=n-r-1,
$
a contradiction to Claim 1.

Thus, $k=2$. As $H$ has a matching of size $3$ and $n=7$,  
$G$ is isomorphic to $K_{4,3}$,
a contradiction to the assumption of $G$.

Claim~\thecountclaim\ holds.
\claimend

\inclaim $r\ne 3$.

\proof
Suppose that $r=3$. By Claim 6, $s=3$ and $|V(H_2)|=|V(H_3)|=3$ by Claims 2 and 3,
implying that $H_i$ is isomorphic to $K_3$ or $P_3$ (the path of order $3$)
for $i=2,3$.

We first show that $|V(H_1)|=3$. 
Suppose that $|V(H_1)|\ge 4$.

For $2\le i\le 3$, 
as $|V(H_1)|>|V(H_i)|=3$, by Lemma~\ref{le6-4}, 
there exists a vertex $u_i\in V(H_i)$ such that 
$|N_G(u_i)\cap V(H_1)|\ge 
1+|N_G(u_i)\cap V(H_i)|\ge |E(H_i)|$.
Let $H_0=G[V(H_1)\cup \{u_2,u_3\}]$.
Observe that 
$$
|E(H_0)|=|E(H_1)|+\sum_{i=2}^3 |N_G(u_i)\cap V(H_1)|
\ge |E(H_1)|+|E(H_2)|+|E(H_3)|=|E(H)|
$$
and
$|V(H_0)|=|V(H_1)|+2=n-4=n-r-1$, 
a contradiction to Claim 1.

Hence $|V(H_1)|=3$, implying that $n=9$ and $6\le |E(H)|\le 9$.

If $k\ge 4$, then $G$ must have a subgraph $H_0$ isomorphic to 
$K_5-e$, the graph obtained from $K_5$ by removing one edge.
Note that $|V(H_0)|=5=n-4=n-r-1$ and $|E(H_0)|=9\ge |E(H)|$,
a contradiction to Claim 1.

Thus, $k\le 3$. If $p_2=1$, then $p_3=1$,
a contradiction to the fact that 
$H$ has three vertex-disjoint paths
of length $2$. 
Hence $p_2\ge 2$.

If $|E(H)|=9$, then $H_i\cong K_3$ 
for all $i=1,2,3$, implying that $k=3$ and 
$G\cong K_{3,3,3}$, 
a contradiction to the assumption of $G$.

Now assume that $|E(H)|\le 8$. 
If $k=3$, as $p_2\ge 2$, 
$G$ has a subgraph $H_0$ isomorphic to $K_{2,2,1}$,
implying that $|E(H_0)|\ge 8\ge |E(H)|$.
But, $|V(H_0)|=5=n-4=n-r-1$, 
a contradiction to Claim 1. 

Hence $k=2$, and each $H_i$ is isomorphic to $P_3$ 
(i.e., the path graph of length $2$),
implying that $|E(H)|=6$.
As $p_2\ge 2$, $G$ has a subgraph $H_0$ isomorphic to $K_{3,2}$,
implying that $|E(H_0)|=6=|E(H)|$
and $|V(H_0)|=5=n-4=n-r-1$, 
a contradiction to Claim 1. 

Claim~\thecountclaim\ holds.
\claimend

As $r=n-q$,
the result follows from Claims 4, 5, 6 and 7.
\proofend

\section{To find $\min\limits_{S\in \P_{r}(G)}
\left |E_G(S)\right |$
for a complete multi-partite graph $G$
\label{sec6-2}
}

In this section, 
we shall show that for a complete multi-partite graph $G$, 
$\min\limits_{S\in \P_{r}(G)}\left |E_G(S)\right |$ 
is equal to $|E_G(S_0)|$ 
for every set $S_0\in \P_r(G)$ 
obtained by a simple algorithm (i.e. Algorithm A).
%for a complete multi-partite graph $G$, we shall provide a simple algorithm to identify a set $S$ in $\P_{r}(G)$ which achieves the minimum value.

For any graph $G$ of order $n$ 
and any integer $r$ with $1\le r\le n$, 
let $\P^*_r(G)$ denote the family of subsets $S\in \P_r(G)$ 
obtained by the following ``greedy" algorithm.
Note that this algorithm chooses 
one vertex by one vertex and always chooses a vertex 
of the minimum degree in the remaining graph.

\noindent {\bf Algorithm A:}
\begin{enumerate}
\item[Step 1] Set $i:=1$, $S:=\emptyset$ and $G_1:=G$;
\item[Step 2] Choose a vertex $u$ from $G_i$ of degree $\delta(G_i)$;
\item[Step 3] Set $S:=S\cup \{u\}$;
\item[Step 4] If $i=r$, then output $S$ and STOP; otherwise, 
set $i:=i+1$ and $G_i:=G_{i-1}\setminus \{u\}$, and return to Step 2.
\end{enumerate}

In the following, we shall show that if $G$ is a complete multi-partite graph,
then $\P^*_r(G)$ is exactly the set of elements $S_0\in \P_r(G)$ such that 
\equ{eq6-14}
{
%\left |\bigcup_{v\in S_0}E_G(v)\right |
|E_G(S_0)|
=\min_{S\in \P_r(G)}
|E_G(S)|.
%\left |\bigcup_{v\in S}E_G(v)\right |.
}

We first establish the following conclusion.

\begin{lemm}\label{le6-5}
Let $G$ be a complete $k$-partite graph of order $n$
with partite sets $X_1,X_2,\cdots,X_k$, where $k\ge 2$.
Assume that $T_0\in \P_b(G)$, where $1\le b\le n$.
Then, 
$|E(G[T_0])|\ge |E(G[T])|$ holds for every $T\in \P_b(G)$
if and only if 
for each pair of distinct numbers $i,j\in \{1,2,\cdots,k\}$, 
if $|X_i\cap T_0|\ge |X_j\cap T_0|+2$, then $X_j\subseteq T_0$.
%if $|X_i\setminus S|\ge 2+|X_j\setminus S|$,
%then $X_j\cap S=\emptyset$.
%either $X_j\subseteq T_0$ or $|X_j\cap T_0|\ge |X_i\cap T_0| -1$.
\end{lemm}

\proof Let $p_i=|X_i|$ for all $i=1,2,\cdots,k$.
Assume that $p_1\ge p_2\ge \cdots \ge p_k\ge 1$.
Let $a_0=0$, and for $j=1,2,\cdots,p_1$,
let $a_j=\max\{1\le i\le k: p_i\ge j\}$.
For example, if $k=3$, $p_1=4, p_2=3$ and $p_3=1$, 
then, $a_1=3, a_2=2, a_3=2$ and $a_4=1$.
In order to complete the proof, we shall show the 
equivalence of the following three statements for any $T_0\in \P_b(G)$:
\begin{enumerate}
\item $|E(G[T_0])|\ge |E(G[T])|$ holds for every $T\in \P_b(G)$;

\item for each pair of distinct numbers $i,j\in \{1,2,\cdots,k\}$, 
if $|X_i\cap T_0|\ge |X_j\cap T_0|+2$, then $X_j\subseteq T_0$.
%if $|X_i\setminus S|\ge 2+|X_j\setminus S|$, then $X_j\cap S=\emptyset$;
%either $X_j\subseteq T_0$ or $|X_j\cap T_0|\ge |X_i\cap T_0| -1$;

\item 
let $h$ be the unique number in $\{1,2,\cdots,p_1\}$ 
determined by the inequality:
$a_1+a_2+\cdots+a_{h-1}<b\le a_1+a_2+\cdots+a_{h-1}+a_h$.
Then, $G[T_0]$ is isomorphic to $K_{q_1,q_2,\cdots,q_k}$,
where $q_i=p_i$ when $p_i\le h-2$, 
and $h-1\le q_i\le h$ otherwise. 
Furthermore, there are exactly 
$b-(a_1+\cdots+a_{h-1})$ indices $i: 1\le i\le k$ 
such that $q_i=h$.
%$|\{1\le i\le k: q_i=h\}|=b-(a_1+\cdots+a_{h-1})$.
\end{enumerate}

(i) $\Rightarrow$ (ii).
Assume that   $T_0\in \P_b(G)$ such that  
$|E(G[T_0])|\ge |E(G[T])|$ holds for every $T\in \P_b(G)$.
As $G$ is a complete $k$-partite graph,
for any $U\subseteq V(G)$,
$G[U]$ is also a complete multi-partite graph, and 
its size is 
\equ{le6-5-e1}
{
|E(G[U])|={|U|\choose 2}
-\sum_{i=1}^k {|U\cap X_i|\choose 2}.
}
Suppose that there is a pair of distinct numbers $i,j\in \{1,2,\cdots,k\}$, 
such that $X_j\not \subseteq T_0$ and $|X_j\cap T_0|\le |X_i\cap T_0| -2$.
Let $v\in X_j\setminus T_0$  and $v'\in X_i\cap T_0$.
For $T'=(T_0\setminus \{v'\})\cup \{v\}\in \P_b(G)$,
by (\ref{le6-5-e1}),  
\eqn{le6-5-e2}
{
|E(G[T'])|-|E(G[T_0])|
&=&{|T_0\cap X_i|\choose 2}+{|T_0\cap X_j|\choose 2}
-{|T_0\cap X_i|-1\choose 2}-{|T_0\cap X_j|+1\choose 2}
\nonumber \\
&=&|T_0\cap X_i|-|T_0\cap X_j|-1>0,
}
a contradiction to the given condition.
Thus, (i) $\Rightarrow$ (ii) holds.

(ii) $\Rightarrow$ (iii).
Without loss of generality, 
assume that $h:=|X_1\cap T_0|\ge |X_j\cap T_0|$
for all $j=2,3,\cdots,k$.
By the given condition, for any $2\le j\le k$,
$|X_j\cap T_0|\ge |X_1\cap T_0| -1=h-1$,
unless $X_j\subseteq T_0$.
Thus, if $p_j\le h-2$, then $X_j\subseteq T_0$.
Obviously, $h$ is the unique number satisfying the inequality below:
$$
a_1+a_2+\cdots+a_{h-1}<b\le a_1+a_2+\cdots+a_{h-1}+a_h,
$$
and
\equ{le6-5-e5}
{
|\{1\le j\le k: |X_j\cap T_0|=h\}|
=b-(a_1+a_2+\cdots+a_{h-1}).
}
Note that $G[T_0]$ is a complete multi-partite graph with partite sets
$T_0\cap X_i$ for $i=1,2,\cdots,k$, 
where $|T_0\cap X_i|=|X_i|=p_i$ whenever $p_i\le h-2$,
and $h-1\le |T_0\cap X_i|\le h$ otherwise. 
Furthermore, by (\ref{le6-5-e5}),
$G[T_0]$ has exactly $b-(a_1+a_2+\cdots+a_{h-1})$ partite sets 
$T_0\cap X_j$ of size $h$.
Thus (iii) holds.

(iii) $\Rightarrow$ (i).
Assume that condition (iii) is satisfied for $T_0$.
Assume that $T'\in \P_b(G)$ such that 
$|E(G[T'])|\ge |E(G[T])|$ holds every $T\in \P_b(G)$.
As (iii) follows from (i), 
condition (iii) is satisfied for $T'$.
As both $T_0$ and $T'$ satisfy condition (iii), 
we have $G[T']\cong G[T_0]$,
implying that $|E(G[T'])|=|E(G[T_0])|$,
and thus 
$|E(G[T_0])|\ge |E(G[T])|$ holds for every $T\in \P_b(G)$.
Hence (iii) $\Rightarrow$ (i) holds.

Therefore (i) $\Leftrightarrow $ (ii) holds, and
the result is proven.
\proofend

\begin{theo}\label{th6-3}
Let $G$ be a complete multi-partite graph of order $n$ with at least two partite sets.
For any $1\le r\le n$ and $S_0\in \P_r(G)$,
$S_0\in \P^*_r(G)$ if and only if 
\equ{th6-3-e0}
{
%\left |\bigcup_{v\in S_0}E_G(u)\right |
|E_G(S_0)|
=
\min_{S\in \P_r(G)}|E_G(S)|.
%\left |\bigcup_{v\in S}E_G(u)\right |.
}
\end{theo}

\proof Let $X_1,X_2,\cdots,X_k$ be the partite sets of $G$, where $k\ge 2.$
%with $|X_i|=p_i$ for all $i=1,2,\cdots,k$.

We need only to show that 
the following statements are equivalent for any $S\in \P_r(G)$:
\begin{enumerate}
\item $S\in \P^*_r(G)$;
\item for each pair of distinct numbers $i,j$ in $\{1,2,\cdots,k\}$, 
if $|X_i\setminus S|\ge 2+|X_j\setminus S|$, then $X_j\cap S=\emptyset$;
%either $X_j\cap S=\emptyset$
%or $|X_i\setminus S|\le 1+|X_j\setminus S|$;
%$|V(G\setminus S_0)\cap X_i|\le 1+|V(G\setminus S_0)\cap X_j|$, 

\item $|E(G\setminus S)|\ge |E(G[T])|$ for each $T\in \P_{n-r}(G)$;
and 

\item 
%$|E_G(S)|\left |\bigcup_{v\in S}E_G(u)\right |\le \left |\bigcup_{v\in S'}E_G(u)\right |$ 
$|E_G(S)|\le |E_G(S')|$
for each $S'\in \P_r(G)$.
\end{enumerate}
The equivalence of (ii) and (iii)
follows from Lemma~\ref{le6-5} by taking $T_0=V(G)\setminus S$,
while (iii) and (iv) are equivalent by (\ref{eq6-08}).
It remains to prove that (i) and (ii) are equivalent.

Assume that $S=\{u_1,u_2,\cdots,u_r\}\in \P^*_r(G)$,
where $u_s$ is the $s$-th vertex in $S$ selected by Algorithm A
for $s=1,2,\cdots,r$, and let $S_s:=\{u_1,u_2,\cdots,u_s\}$
and $G_s:=G\setminus \{u_1,\cdots, u_{s-1}\}$. 

Suppose that (ii) fails. 
Then, there exist distinct numbers $i,j$ in $\{1,2,\cdots,k\}$
such that $X_j\cap S\ne \emptyset$ and 
$|X_i\setminus S|\ge 2+|X_j\setminus S|$.
As $X_j\cap S\ne \emptyset$, there exists $u_q\in S\cap X_j$.
If $|X_j\cap S|\ge 2$, we assume $u_q$ is chosen with the largest possible 
value of $q$.

As $q$ is the largest number in $\{1,2,\cdots,r\}$ 
such that $u_q\in X_j$,
%Observe that 
$$
|X_i\setminus S_{q-1}|\ge |X_i\setminus S|\ge 
2+|X_j\setminus S|=1+|X_j\setminus S_{q-1}|.
$$
Thus, for any $u\in X_i\setminus S_{q-1}$, we have 
	$$
	d_{G_q}(u_q)=|V(G_q)|-|X_j\setminus S_{q-1}|
	\ge |V(G_q)|-|X_i\setminus S_{q-1}|+1
	=d_{G_q}(u)+1,
	$$
a contradiction to the condition that 
$u_q$ has the minimum degree in $G_q$.
Thus (i) $\Rightarrow$ (ii) holds.

Assume that condition (ii) is satisfied.
% for each pair of distinct numbers $i,j$ in $\{1,2,\cdots,k\}$, 
%if $|X_i\setminus S|\ge 2+|X_j\setminus S|$, then $X_j\cap S=\emptyset$;
Then for each pair of distinct numbers $j,j'$ in $\{1,2,\cdots,k\}$ such that 
$X_j\cap S\ne \emptyset$ and $X_{j'}\cap S\ne \emptyset$, the difference between
$|X_j\setminus S|$ and $|X_{j'}\setminus S|$ is at most one. Moreover, if $X_i\cap S= \emptyset$ for some $i\in \{1,2,\cdots,k\}$, then 
$$|X_i|=|X_i\setminus S|\le 1+\min\{|X_j\setminus S|:X_j\cap S\ne \emptyset, 1\le j \le k\}.$$

Note that the vertices in $S$ can be determined by the 
following Algorithm.

\noindent {\bf Algorithm B}:
\begin{enumerate}
	\item[Step 1.] Set $t:=r$ and $S_t:=S$; % $G_{r-1}:=G\setminus S$;
	\item[Step 2.] 
choose $j:1\le j\le k$ such that 
	$X_j\cap S_t\ne \emptyset$ and 
	$|X_j\setminus S_t|\le |X_{j'}\setminus S_t|$
	for each $j'\ne j$ with $X_{j'}\cap S_t\ne \emptyset$;
	\iffalse 
	if 
	$\sum\limits_{s=2}^k |S_t\cap X_s|(|X_1\setminus S_t|-|X_s\setminus S_t|) =0$, 
	let $j=1$;
	otherwise,
	choose a number $j:2\le j\le k$ such that 
	$S_t\cap X_j\ne \emptyset$ and 
	$|X_1\setminus S_t|>|X_j\setminus S_t|$. % has the maximum value;
	\fi
	
	\item[Step 3.] let $u_t$ be a vertex in $S_t\cap X_j$;
	
	\item[Step 4.] if $t=1$, then STOP;
	otherwise, set $t:=t-1$, $S_t:=S_{t+1}\setminus \{u_{t+1}\}$,
	and go to Step 2.
\end{enumerate}

Assume that $u_1,u_2,\cdots,u_r$ are vertices in $S$ determined 
by Algorithm B.
Note that for each $t=1,2,\cdots,r$,
$S_t=\{u_1,\cdots,u_t\}$ and 
$u_t$ has the minimum degree in $G\setminus S_{t-1}$,
where $S_0=\emptyset$.
Hence, $S\in \P^*_r(G)$
and (ii) $\Rightarrow$ (i) holds.

The result is proven.
\proofend

\begin{coro}\label{co6-4}
Let $G$ denote $K_{p_1,p_2,\cdots,p_k}$
and $n=p_1+p_2+\cdots+p_k$,
where $k\ge 2$.
For $2\le q\le n-1$, if 
$
\ell_q(G)=|E(G)|-\min_{S\in \P_{n-q+1}(G)}
\left |E_G(S)\right |,
$
then,  
\equ{cp6-4-e1}
{
ar(G,\T_q)=1+\sum_{1\le i<j\le k}p_ip_j
-|E_G(S_0)|,
%\left |\bigcup_{v\in S}E_G(v)\right |,
}
where $S_0\in \P^*_{n-q+1}(G)$.
% and $|E_G(S_0)|$ is givenin (\ref{eq6-20}), (\ref{eq6-21}) and (\ref{eq6-22})respectively for $r=1,2$ and $3$.
\end{coro}

\proof 
The result follows from Theorems~\ref{theoq} and~\ref{th6-3}.
\proofend

\section{To determine $ar(K_{p_1,\cdots,p_k},\T_q)$ for $n-3\le q\le n-1$
\label{sec6-3}
}

In this section, 
we shall 
%apply Theorem~\ref{th6-2} and
%Corollary~\ref{co6-4} to 
give an explicit expression for 
$ar(K_{p_1,\cdots,p_k},\T_q)$ 
when $\max\{2,n-3\}\le q\le n-1$, where $n=p_1+\cdots+p_k$.

If $(G,q)$ is an ordered pair in (\ref{eq6-13}),
$ar(G,\T_q)$ can be determined by applying 
Theorems~\ref{theoq} and~\ref{co6-1} (ii) that 
%It can be verified that 
\equ{bth6-3}
{
ar(G,\T_q)=1+\ell_q(G)=
\left \{
\begin{array}{ll}
4, \qquad &\mbox{if }(G,q)=(K_{3,3},4) \mbox{ or }(K_{4,3},4);\\
%4,  &\mbox{if }(G,q)=(K_{4,3},4);\\
10,  &\mbox{if }(G,q)=(K_{3,3,3},6).
\end{array}
\right.
}
For example, if $(G,q)=(K_{4,3},4)$
and $H$ is a member in $\msq{K_{4,3}}{4}$ such that 
$\seq(H)\succeq \seq(H')$ for every $H'\in \msq{K_{4,3}}{4}$,
then, by Theorem~\ref{co6-1} (ii),
$\seq(H)$ is either $(2,2,2,1)$  
or $(3,1,1,1,1)$,
implying that $\ellq{K_{4,3}}{4}=3$.
Similarly, if $(G,q)=(K_{3,3,3},6)$
and $H$ is a member in $\msq{K_{3,3,3}}{6}$ such that 
$\seq(H)\succeq \seq(H')$ for every $H'\in \msq{K_{3,3,3}}{6}$,
then, by Theorem~\ref{co6-1} (ii),
$\seq(H)$ is one of the sequences 
$(3,3,3)$, $(4,2,2,1)$ or $(5,1,1,1,1)$,
implying that $\ellq{K_{3,3,3}}{6}=9$.

In the following, we consider the case that 
$G$ is a complete multi-partite graph of order $n$ 
and $\max\{2,n-3\}\le q\le n-1$ such that 
$(G,q)$ is not an ordered pair in (\ref{eq6-13}).

Note that for any complete multi-partite graph $G$ and 
$u\in V(G)$, $d_G(u)=\delta(G)$ if and only if 
$u$ is contained in a partite set with the largest cardinality. 
Now let $G=K_{p_1,p_2,\cdots,p_k}$
with partite sets $X_1,X_2,\cdots,X_k$
and $n=|V(G)|$,
where $k\ge 2$, $p_1\ge p_2\ge \cdots \ge p_k\ge 1$
and $|X_i|=p_i$ for all $i=1,2,\cdots,k$.
%Clearly,a vertex $u$ in $G$ has the minimum degree in $G$ if and only if  $u\in X_j$ such that $|X_j|=\max\limits_{1\le i\le k}|X_i|$.
Clearly, for any $S=\{x_1,x_2,\cdots,x_s\}\subseteq V(G)$, 
if $x_i\in X_{j_i}$ for all $i\in [s]$, then 
\equ{eq6-19} 
{
|E_G(S)|=\sum_{i=1}^s (n-p_{j_i})-\sum_{1\le i_1<i_2\le s}\sigma (j_{i_2}-j_{i_1})
=sn-\sum_{i=1}^s p_{j_i}-\sum_{1\le i_1<i_2\le s}\sigma (j_{i_2}-j_{i_1}),
}
where $\sigma(x)$ is the function defined by $\sigma(0)=0$ 
and $\sigma(x)=1$ when $x\ne 0$.
Then, for $S\in \P^*_t(G)$, where $2\le t\le 4$,
%By applying Algorithm A and (\ref{eq6-19}), 
%$\left |\bigcup_{v\in S}E_G(v)\right |$ 
$|E_G(S)|$ can be determined by applying Algorithm A and (\ref{eq6-19}) 
as follows.
If $S\in \P^*_2(G)$,
\equ{eq6-20}
{
|E_G(S)|
=
\left \{
\begin{array}{ll}
2n-2p_1,        \quad &\mbox{if }p_1>p_2;\\
2n-p_1-p_2-1, \quad  &\mbox{if }p_1=p_2.
\end{array}
\right.
}
If $S\in \P^*_3(G)$,
\equ{eq6-21}
{
%\left |\bigcup_{v\in S}E_G(v)\right |
|E_G(S)|
=
\left \{
\begin{array}{ll}
3n-3p_1,                  \quad &\mbox{if }p_1\ge p_2+2;\\
3n-2p_1-p_2-2,        \quad &\mbox{if }p_2+1\ge p_1\ge p_3+1;\\
%p_1=p_2+1\mbox{ or }p_1=p_2\ge p_3+1;\\
3n-p_1-p_2-p_3-3, \quad &\mbox{if }p_1=p_2=p_3.
\end{array}
\right.
}
If $S\in \P^*_4(G)$,
	\equ{eq6-22}
	{
		%\left |\bigcup_{v\in S}E_G(v)\right |
		|E_G(S)|
		=
		\left \{
		\begin{array}{ll}
			4n-4p_1,   \qquad &\mbox{if } p_1\ge p_2+3;\\
			4n-3p_1-p_2-3, \quad 
&\mbox{if } p_1=p_2+2 \mbox{ or }p_1=p_2+1\ge p_3+2;\\
%&p_3+2\le p_2+1\le p_1\le p_2+2;\\
4n-2p_1-2p_2-4, \quad &\mbox{if } p_1=p_2\ge p_3+1;\\
4n-2p_1-p_2-p_3-5, 
%&p_2+1\ge p_1\ge p_4+1 \mbox{ and } p_2=p_3;\\
&\mbox{if } p_4+1\le p_1\le p_2+1=p_3+1;\\
%\quad &p_1=p_2+1=p_3+1 \mbox{ or }p_1=p_2=p_3\ge p_4+1;\\
			4n-\sum\limits_{i=1}^4p_i-6,%+(n-p_2)+(n-p_3)+(n-p_4)-6, 
\quad &\mbox{if }p_1=p_2=p_3=p_4.\\
		\end{array}
		\right.
	}
	
\begin{theo}\label{th6-5}
Let $G$ denote $K_{p_1,p_2,\cdots,p_k}$
and $n=p_1+p_2+\cdots+p_k$,
where $k\ge 2$ and $p_1\ge p_2\ge \cdots\ge p_k\ge 1$.
For $\max\{2,n-3\}\le q\le n-1$, 
if $(G,q)$ is an ordered pair in (\ref{eq6-13}),
then $ar(G,\T_{q})$ is given in (\ref{bth6-3}); 
otherwise, 
%if $q\ge 2$ and $(G,q)$ is not an ordered pair in (\ref{eq6-13}),then, 
\equ{th6-3-e1}
{
ar(G,\T_{q})=1+\sum_{1\le i<j\le k}p_ip_j
-|E_G(S_0)|,
%\left |\bigcup_{v\in S}E_G(v)\right |,
}
where $S_0\in \P^*_{n-q+1}(G)$ and $|E_G(S_0)|$ is given
in (\ref{eq6-20}), (\ref{eq6-21}) and (\ref{eq6-22})
respectively for $q=n-1, n-2$ and $n-3$.
%$r=1,2$ and $3$.
\end{theo}

\proof
As $G=K_{p_1,p_2,\cdots,p_k}$
and $(G,q)$ is not an ordered pair in (\ref{eq6-13}),
where $\max\{2,n-3\}\le q\le n-1$,
by Theorem~\ref{th6-2} and Corollary~\ref{co6-4},
  we have 
\equ{th6-3-e2}
{
ar(G,\T_{q})
%=1+\ell_{n-r}(G)=1+E(G)-\min_{S\in \P_{r+1}(G)}|E(G[S])|
=1+E(G)-|E(G[S_0])|,
}
where $S_0\in \P^*_{n-q+1}(G)$.
Thus, the result follows from (\ref{eq6-20}), (\ref{eq6-21}) and (\ref{eq6-22}).
\proofend

\section{To determine $ar(K_{p_1,\cdots,p_k},\T_q)$ for $(4n-2)/5\le q\le n-1$
\label{sec6-4}}

Let $G$ denote $K_{p_1,\cdots,p_k}$.
In this section, we consider the case 
when $(4n-2)/5\le q\le n-1$,
where $n=p_1+p_2+\cdots+p_k$. 
We shall show that in this case,  
$ar(G,\T_q)$
can be calculated by $|E(G)|-|E_G(S_0)|+1$ 
for any $S_0\in \P^*_{n-q+1}(G)$
in Theorem~\ref{th6-6}.
Moreover, we give the explicit expression 
for $ar(K_{p_1,\cdots,p_k},\T_q)$ in Corollary~\ref{th6-6-co1} with a further condition that $p_1-p_2\ge (n+2)/5$.

\begin{lemm}\label{le6-6}
Let $H$ be a spanning subgraph of a complete multi-partite 
graph $G$ of order $n$
and let $r$ be a positive integer such that $n\ge 5r-2$.
%If $H$ has $s$ components $H_1,H_2,\cdots, H_s$
%with $|V(H_1)|\ge \cdots \ge |V(H_s)|$, where $s\ge 3$, 
%such that 
If $com(H)\ge 2$, $or_2(H)\le r$
and $n-(or_1(H)+or_2(H))\le r$,
then there exists $S\subseteq V(H)\setminus V(H_1)$ 
with $|S|\le or_2(H)-1$ 
such that 
$
|E(G[V(H_1)\cup S])|\ge |E(H)|,
$
where $H_1$ is a largest component of $H$.
\end{lemm}

\proof %Let $h=or_2(H)$.
We prove this result by induction on $or_2(H)$.
It is trivial if $or_2(H)=1$.
Assume that it holds when $or_2(H)<d$, where $d\ge 2$.
Now consider the case that $or_2(H)=d$.

Let $H_2,H_3,\cdots,H_s$ be the components of $H$
different from $H_1$
and
let $W$ denote the subgraph $H\setminus V(H_1)$.
Clearly, $W$ consists of components $H_2,H_3,\cdots,H_s$
and $|V(W)|=|V(H_2)|+(|V(H_3)|+\cdots+|V(H_s)|)
\le or_2(H)+r\le 2r$, implying that 
$|V(H_1)|\ge n-2r\ge 3r-2$.

\inclaim
There exist $u$ and $U$ with $u\in U\subseteq V(W)$ 
%and $U\subseteq V(W)$ with $u\in W$ 
such that $|U\cap V(H_i)|=1$ for each $i\in [s]\setminus \{1\}$
and
$$
|N_G(u)\cap V(H_1)|\ge \sum_{v\in U} d_W(v).
%\frac{t-1}{t}\times |V(H_1)|,
$$
%where $t=\chi(W)$.

\proof %Let $V_1=V(H_1)$ and $V_2=V(W)$. 
%where $W'$ consists of non-trivial components of $W$.
Let $t=\chi(W)$. If $t=1$, then $W$ is an 
empty graph and Claim 1 is trivial. 
Now assume that $t\ge 2$.
By Lemma~\ref{le6-40}, there exists $u\in V(W)$
such that  
\equ{le6-6-e1}
{
|N_G(u)\cap V(H_1)|\ge \frac{t-1}{t}\times |V(H_1)|.
}

%Without loss of generality,
Assume that $u\in V(H_b)$, where $2\le b\le s$. 
Clearly, $d_W(u)\le |V(H_b)|-1$.
By Lemma~\ref{le6-3}, there exists $u_i\in V(H_i)$ 
for each $i\in [s]\setminus \{1,b\}$ 
such that 
\equ{le6-6-e3}
{
\sum_{2\le i\le s\atop i\ne b}d_W(u_i)
\le \frac {t'-1}{t'} 
|V(W)\setminus V(H_b)|
\le \frac {t-1}{t} (|V(W)|-|V(H_b)|),
}
where $t'=\chi(W\setminus V(H_b))$.
Let $U=\{u\}\cup \{u_i: 2\le i\le s, i\ne b\}$.
Since $|V(H_b)|\le r$ and $|V(W)|\le 2r$, we have 
\eqn{le6-6-e2}
{
\sum_{v\in U}d_W(v)
&\le &
|V(H_b)|-1+\frac {t-1}{t} (|V(W)|-|V(H_b)|)
\nonumber \\ 
&\le &
\frac {t-1}{t} |V(W)|+\frac {1}{t} |V(H_b)|-1
\nonumber \\ 
&\le &
\frac {t-1}{t} (2r)+\frac {1}{t} r-1
\nonumber \\
&=&\frac {r(2t-1)}{t}-1.
}
Since $|V(H_1)|\ge 3r-2$,  
\equ{le6-6-e4}
{
\frac {r(2t-1)}{t}-1
=\frac{t-1}{t}(3r-2)-\frac{(t - 2)(r-1)}t 
\le \frac{t-1}{t}(3r-2)
\le 
\frac{t-1}{t}|V(H_1)|.
}
Thus, Claim 1 follows from (\ref{le6-6-e1}),
 (\ref{le6-6-e2})
and (\ref{le6-6-e4}).
\claimend

By Claim 1, there exist $u$ and $U$ with 
$u\in U\subseteq V(W)$ such that the
properties in Claim 1 hold.
Let $H'$ be the spanning subgraph of $G$ 
obtained from $H$  by 
adding all edges in 
$\{uv:v\in N_G(u)\cap V(H_1)\}$
%$E_G(u,V(H_1))$
and deleting all edges in 
$\bigcup_{v\in U}E_W(v)$.

As $d\ge 2$, 
it is clear that $or_2(H')= or_2(H)-1\le r-1$ and 
$or_3(H')+\cdots+or_{com(H')}(H')
= or_3(H)+\cdots+or_s(H)\le r$.

As $|V(H_1')|\ge or_2(H)\ge 2$, $H'_1:=G[V(H_1)\cup \{u\}]$ is connected and  
the largest component of $H'$.
Since 
$or_2(H')=or_2(H)-1=d-1$, 
by inductive assumption, 
there exists $S'\subseteq V(H')\setminus V(H'_1)$
with $|S'|\le or_2(H')-1=or_2(H)-2$ 
such that $|E(G[V(H'_1)\cup S'])|\ge |E(H')|$.

Let $S=S'\cup \{u\}$. 
Then $|S|\le or_2(H)-1$ and 
$V(H'_1)\cup S'
=V(H_1)\cup S$.
By Claim 1, we have 
$|E(H')|\ge |E(H)|$.
Thus, 
$$
|E(G[V(H_1)\cup S])|= |E(G[V(H'_1)\cup S'])|
\ge |E(H')|\ge |E(H)|.
$$
The result follows.
\proofend

\begin{theo}\label{th6-6}
Let $G$ be a complete multi-partite graph of order $n$ with at least two partite sets.
If $q$ is an integer with 
$(4n-2)/5 \le q\le n-1$, then 
\equ{eq6-16}
{
ar(G,\T_q)= |E(G)|+1- |E_G(S_0)|,
%\ell_q(G)= |E(G)|- |E_G(S_0)|.
}
where $S_0$ is a member in $\P^*_{n-q+1}(G)$.
\end{theo}

\proof Let $H\in \msq{G}{q}$ such that 
$or_1(H)\ge or_1(H')$
%$(H) \succeq \seq(H')$ 
holds for every 
$H'\in \msq{G}{q}$.

\inclaim $or_2(H)=1$.

\proof
Suppose that $or_2(H)\ge 2$.
Let $r=n-q$. As $(4n-2)/5 \le q$, we have $n\ge 5r-2$.
Let $s=com(H)$, then $s\ge 3$ by Lemma~\ref{le2-1}.

By Theorem~\ref{co6-1} (i),
$or_1(H)+or_2(H)=q$ and $or_2(H)=or_3(H)$.
%$or_2(H)=or_3(H)=\cdots=or_{s-1}(H)\ge or_s(H)$.
Clearly, 
$$
or_3(H)+\cdots+or_s(H)=n-(or_1(H)+or_2(H))=n-q=r
$$ 
and $2\le or_2(H)=or_3(H)\le r$.

Let $H_1$ be the largest component of $H$.
Since $n\ge 5r-2$, by Lemma~\ref{le6-6},
there exists $S\subseteq V(H)\setminus V(H_1)$ 
with $|S|\le or_2(H)-1$ 
such that 
$|E(G[V(H_1)\cup S])|\ge |E(H)|$.
%Thus, $G[V(H_1)\cup S]$ is connected.
Since $or_2(H)\ge 2$, $|E(G[V(H_1)\cup S])|\ge |E(H)|$
implies that $|S|\ge 1$.
Note that $G[V(H_1)\cup S]$ is connected as $|V(H_1)|\ge or_2(H)\ge 2$, and
$$
|V(H_1)|+|S|\le or_1(H)+or_2(H)-1=q-1,
$$
by which the spanning subgraph $H'$ of $G$ 
with edge set $E(G[V(H_1)\cup S])$ 
belongs to $\sq{G}{q}$.
As $|E(G[V(H_1)\cup S])|\ge |E(H)|$, 
$H'\in \msq{G}{q}$, while $or_1(H')=|V(H_1)|+|S|>|V(H_1)|$,
a contradiction to the assumption that 
$or_1(H)\ge or_1(H')$.
Hence $or_2(H)=1$.
\claimend

By Theorem~\ref{co6-1} (i) and Claim 1,  $|V(H_1)|=q-or_2(H)=q-1$
and 
$$
\ell_q(G)=|E(H)|=|E(H_1)|=|E(G)|-|E_G(S_0)|
$$
where $S_0=V(H)\setminus V(H_1)\in \P_{n-q+1}(G)$.
By (\ref{eq6-8}), $|E_G(S_0)|\le |E_G(S)|$ 
for each $S\in \P_{n-q+1}(G)$,
thus by Theorem~\ref{th6-3}, $S_0\in \P^*_{n-q+1}(G)$.

Due to Theorem~\ref{theoq}, $ar(G,\T_q)=\ellq{G}{q}+1$.
The result then follows.
\proofend

\begin{coro}\label{th6-6-co1}
	Let $G$ denote $K_{p_1,p_2,\cdots,p_k}$
	and $n=p_1+p_2+\cdots+p_k$,
	where $k\ge 2$ and $p_1\ge p_2\ge \cdots\ge p_k\ge 1$.
If $(4n-2)/5\le q\le n-1$ and 
$p_1-p_2\ge (n+2)/5$, then %n+1-(4n-2)/5$, then 
$$
ar(G,\T_q)=|E(G)|+1-(n-q+1)(n-p_1).
$$ 
\end{coro}

\proof By Theorem~\ref{th6-6},
$ar(G,\T_q)=1+|E(G)|-|E_G(S)|$
for any $S\in \P^*_{n-q+1}(G)$.
Since $p_1-p_2\ge (n+2)/5 = n-(4n-2)/5\ge n-q$,
by Algorithm A, 
there exists $S_0\in \P^*_{n-q+1}(G)$
with $S_0\subseteq X_1$,
where $X_1$ is the partite set of $G$ 
with $|X_1|=p_1$. 

Since $S_0\subseteq X_1$ and $|S_0|=n-q+1$, we have $|E_G(S_0)|=(n-p_1)(n-q+1)$.
Thus, the result holds.
\proofend

\noindent \textbf{Remark.}
Note that the ordered pairs $(K_{3,3}, 4)$ and $(K_{3,3,3}, 6)$ in (\ref{eq6-13}) imply that the condition $q\ge (4n-2)/5$ in Theorem~\ref{th6-6} cannot be improved into $q\ge \lfloor (4n-2)/5\rfloor$,
but we guess it
%and Corollary~\ref{th6-6-co1} 
may be true when $q\ge (2n+1)/3$.
In addition, when $q\le n-4$ and $q$ is close to $n$,
the value of $ar(K_{p_1,\cdots,p_k},\T_{q})$ might also be obtained
by applying Corollary~\ref{co6-4}. 
However, the exact values of $ar(K_{p_1,\cdots,p_k},\T_{q})$ for all $q$ are generally hard to compute, even for complete bipartite graphs, 
as also mentioned in~\cite{5}.

\end{document}